\documentclass[a4paper,12pt]{article}

\usepackage{graphicx}
\RequirePackage{amsmath}
\RequirePackage{amssymb}

\renewcommand\({\left(}
\renewcommand\){\right)}

\newcommand{\half}{{\mbox{$\frac{1}{2}$}}}

\renewcommand{\[}{\begin{equation}}
\renewcommand{\]}{\end{equation}}

\newcommand{\OU}{Ornstein--Uhlenbeck}

\newlength{\figurewidth}
\let\epsilon=\varepsilon

\renewcommand\baselinestretch2

\begin{document}

   \title{Stochastic Stokes' drift of a flexible dumbbell}

   \author{Kalvis M. Jansons\footnote{Email:
       \texttt{stokesdrift@kalvis.com}} \\ Department of Mathematics,
     University College London, \\ Gower Street, London WC1E 6BT, UK}
   \label{firstpage}
   
   
   \maketitle

   \begin{abstract}
     We consider the stochastic Stokes' drift of a flexible dumbbell.
     The dumbbell consists of two isotropic Brownian particles
     connected by a linear spring with zero natural length, and is
     advected by a sinusoidal wave.  We find an asymptotic
     approximation for the Stokes' drift in the limit of a weak wave,
     and find good agreement with the results of a Monte Carlo
     simulation.  We show that it is possible to use this effect to
     sort particles by their flexibility even when all the particles
     have the same diffusivity.
   \end{abstract}

   \section{Introduction}
   \label{intro}
   
   Stochastic Stokes' drift is the modification of classical Stokes'
   drift\cite{Stokes} by random forcing, for example Brownian
   motion\cite{Jansons,Jansons2}.  This process can be used for
   particle sorting\cite{Jansons,Jansons2} as, to leading order in the
   strength of the wave motion, the drift from waves with different
   spatial or temporal frequencies are additive.  This is because the
   Stokes' drift first appears at second order in the wave strength
   and the cross-terms average to zero. Thus it is possible to make
   particles of different diffusivities drift in different directions.
   There are also oceanographic applications\cite{Restrepo}. Some
   exact results are known in one-dimension\cite{Broeck, Bena}.
   
   Processes of this kind, however, are not restricted to the direct
   fluid mechanical applications, and occur in many other systems, for
   example in molecular motors, see \cite{Reimann} and the references
   therein.
   
   In this study, we consider a simple model of a flexible particle,
   namely an elastic dumbbell, for which analytical results are still
   possible.  Such a dumbbell model might be used as a crude model of
   a flexible polymer in solution.  We show that particles of
   different flexibilities (i.e. different spring constants) have
   different Stokes' drift velocities even if they all have the same
   diffusivity.  

   Furthermore, there exist parameters for which the dependence of the
   Stokes' drift on the flexibility is non-monotonic, even for
   particles having the same diffusivities.  So, even in one-dimension,
   particle selection is possible.

   Because, to leading order in the Stokes' drift, waves of different
   spatial or temporal frequencies do not interact, the different
   drifts for particles having different flexibilities but the same
   diffusivities can be made to `fanout' (i.e. drift in different
   directions). This extends the possibilities for using stochastic
   Stokes' drift in particle sorting.

   \section{The stochastic dumbbell process}
   \label{dumbbell}
   
   Consider the motion of two identical isotropic Brownian point
   particles $X$ and $Y$, with no hydrodynamic interactions between them,
   and joined by a linear spring, with relaxation time $\lambda^{-1}$,
   and of zero natural length. We henceforth refer to $\lambda$ as the
   spring constant (even though it is not in its usual units). 

   The absence of hydrodynamic interactions would not be expected to
   qualitatively change the results, and including them would lead to
   a very much harder mathematical problem.  Also, if hydrodynamic
   interactions were included, one would probably want to consider a
   more realistic model of a flexible polymer than the dumbbell model
   considered here. This, however, is left for future work.

   In the absence of a wave, i.e.\ in a fluid at rest, we set the
   diffusivity of the centre of the two particles to $D \equiv
   \half\sigma^2$.  Thus, the stochastic dumbbell equations, in It\^o
   form, are
   \[\label{dumbbell1}
   \begin{aligned}
   dX_t &= \(\epsilon f(X_t, t) - \half \lambda (X_t - Y_t)\)dt 
   + 2^{1/2} \sigma dB_t,\\
   dY_t &= \(\epsilon f(Y_t, t) + \half \lambda (X_t - Y_t)\)dt 
   + 2^{1/2} \sigma dW_t,
   \end{aligned}
   \]
   where $B$ and $W$ are independent standard Brownian motions, $f$
   represents the wave motion, and $\epsilon$ is a dimensionless
   parameter, which will be used as a small parameter in the
   asymptotic analysis to follow.

   This system can be studied in several different ways.  We adopt a
   rational asymptotic approach, as often used in fluid mechanics,
   albeit using stochastic calculus, rather than, say, centre manifold
   theory. The structure of the calculation is simple, as it uses 
   well-known results throughout, and is very similar to the standard
   argument for classical Stokes' drift.
   
   Suppose that the time $t \in [-T, T]$, where
   $T=\epsilon^{-\frac12}T_0$, and $T_0$ is the largest natural time
   scale in the system.  These natural time scales are the spring
   relaxation time $\lambda^{-1}$, the period of the wave forcing, and
   the typical time for a particle to diffuse a wavelength.

   Now consider two new processes given by $U_t \equiv \half (X_t -
   Y_t)$, and $V_t \equiv \half (X_t + Y_t)$, where `$\equiv$' means
   defined to be equal.  From \eqref{dumbbell1}, we find
   \[\label{dumbbell2}
   \begin{aligned}
   dU_t &= \half\epsilon (f(X_t, t) - f(Y_t, t)) dt
   - \lambda U_t dt + \sigma d\bar{B}_t, \\
   dV_t &= \half\epsilon (f(X_t, t) + f(Y_t, t)) dt
   + \sigma d\bar{W}_t,
   \end{aligned}
   \]
   where $\bar B$ and $\bar W$ are independent standard Brownian
   motions.

   We now consider a formal asymptotic expansion in $\epsilon$ for
   \[
   \begin{aligned}
     X_t &= X^{(0)}_t + \epsilon X^{(1)}_t + \epsilon^2 X^{(2)}_t +
     \cdots, \qquad
     Y_t &= Y^{(0)}_t + \epsilon Y^{(1)}_t + \epsilon^2 Y^{(2)}_t +
     \cdots, \\
     U_t &= U^{(0)}_t + \epsilon U^{(1)}_t + \epsilon^2 U^{(2)}_t +
     \cdots, \qquad
     V_t &= V^{(0)}_t + \epsilon V^{(1)}_t + \epsilon^2 V^{(2)}_t +
     \cdots, \\
   \end{aligned}
   \]
   in the limit $\epsilon\to0$.  We now take each order in $\epsilon$
   one at a time.
   
   \subsection{Order $\epsilon^0$}

   At this order, \eqref{dumbbell2} becomes
   \[\label{order0}
   dU^{(0)}_t = -\lambda U^{(0)}_t dt + \sigma d\bar{B}_t, \qquad
   dV^{(0)}_t = \sigma d\bar{W}_t.
   \]
   This implies
   \[\label{sol0}
      U^{(0)}_t = e^{-\lambda t} 
      \bar{\bar B}(\half\lambda^{-1}\sigma^2 e^{2\lambda t}), \qquad
      V^{(0)}_t = \sigma \bar{W}_t,
   \]
   where $\bar{\bar B}$ is a standard Brownian motion, with $\bar{\bar
     B}(0) = 0$, independent of $\bar W$.  Note that, for simplicity,
   we have started the stable \OU{} process $U^{(0)}$ with its
   stationary law; though this choice does not affect the long-time
   behaviour of interest, it does avoid starting transients. So to
   leading order, the centre of the dumbbell moves as a Brownian
   motion with diffusivity $D \equiv \half\sigma^2$, and explains the
   strange choice constants in \eqref{dumbbell1}.
   
   Dropping the bars, we find the following solution of
   \eqref{dumbbell1} to order $\epsilon^0$:
   \[\label{XYsol0}
    \begin{aligned}
      X^{(0)}_t &= \sigma W_t
      + e^{-\lambda t}
      B(\half\lambda^{-1}\sigma^2 e^{2\lambda t}), \\
      Y^{(0)}_t &= \sigma W_t
      - e^{-\lambda t}
      B(\half\lambda^{-1}\sigma^2 e^{2\lambda t}).
    \end{aligned}
   \]
   Note, however, that $B$ and $W$ are not the same Brownian motions
   as in equation \eqref{dumbbell1}.

   \subsection{Order $\epsilon^1$}

   From \eqref{dumbbell2}, at order $\epsilon^1$, we find
   \[\label{order1}
   \begin{aligned}
   \frac{dU^{(1)}_t}{dt}
   &= -\lambda U^{(1)} + \half \(f(X^{(0)}_t, t) - f(Y^{(0)}_t, t)\), \\
   \frac{dV^{(1)}_t}{dt} &= \half \(f(X^{(0)}_t, t) + f(Y^{(0)}_t, t)\),
   \end{aligned}
   \]
   from which the corresponding expressions for $X^{(1)}$ and
   $Y^{(1)}$ immediately follow.  We, however, need only $X^{(1)}$,
   which is given by
   \[\label{x1}
   \frac{dX^{(1)}_t}{dt} = -\lambda U^{(1)}_t + f(X^{(0)}_t, t),
   \]
   where
   \[\label{u1}
   U^{(1)}_t = \half \int_0^\infty e^{-\lambda \alpha}
   \(f(X^{(0)}_{t-\alpha}, t-\alpha) - f(Y^{(0)}_{t-\alpha}, t-\alpha)\)
   d\alpha,
   \]
   where taking the upper limit of the integral to be $\infty$ is
   equivalent to ignoring starting transients.  Even though strictly
   the asymptotic series is valid only on a finite interval $[-T,T]$,
   extending the integral beyond $T$ introduces only exponentially
   small terms that do not contribute to the asymptotic series at any
   finite order in $\epsilon$.  Also, we are interested only in the
   Stokes' drift to leading order, and this appears at order
   $\epsilon^2$, as it does in the classical Stokes' drift problem.

   \subsection{Order $\epsilon^2$}

   From \eqref{dumbbell2}, at order $\epsilon^2$, we find
   \[\label{order2}
   \begin{aligned}
   \frac{dU^{(2)}_t}{dt}
   &= -\lambda U^{(2)} + \half \(f'(X^{(0)}_t, t)X^{(1)}_t 
   - f'(Y^{(0)}_t, t)Y^{(1)}_t\), \\
   \frac{dV^{(2)}_t}{dt} &= \half \(f'(X^{(0)}_t, t)X^{(1)}_t 
   + f'(Y^{(0)}_t, t)Y^{(1)}_t\),
   \end{aligned}
   \]
   where $f'$ is the derivative of $f$ with respect to its first
   argument.  From \eqref{order2}, the corresponding expressions for
   $X^{(2)}$ and $Y^{(2)}$ immediately follow.

   \section{The stochastic Stokes' drift of a dumbbell}
   \label{stokesdrift}

   We now consider a special case in which
   \[
   f(x, t) = u \cos(k x - \omega t + \phi),
   \]
   where $u$, $k$ and $\omega$ are constants, and $\phi$ is a random
   phase, uniformly distributed on $[0, 2\pi)$.  The Stokes' drift will
   not depend on $\phi$, but this choice of phase avoids starting
   transients.
   
   In this section, we are interested only in finding Stokes' drift
   \[\label{drift1}
   {\cal V} \equiv \lim_{t\to\infty} t^{-1} \(V_t - V_0\).
   \]
   From \eqref{order0} and \eqref{order1}, we see that there is no
   contribution to the Stokes' drift at orders $\epsilon^0$ and
   $\epsilon^1$, so we find
   \[\label{LOD}
   {\cal V} = \epsilon^2 E\left[\frac{dV^{(2)}_t}{dt}\right] + \cdots,
   \]
   where we have used the Ergodic Theorem to replace the time limit in
   \eqref{drift1} by an expectation, and this is why it was important
   to avoid starting transients.  Using \eqref{order2}, reduces
   \eqref{LOD} to
   \[\label{drift2}
   {\cal V} = \epsilon^2 E\left[ \half \(f'(X^{(0)}_t, t)X^{(1)}_t
   + f'(Y^{(0)}_t, t)Y^{(1)}_t\)\right] + \cdots.
   \]
   Note that both of the terms in the expectation of \eqref{drift2}
   have the same law, thus
   \[\label{drift3}
   {\cal V} 
   = \epsilon^2 E\left[f'(X^{(0)}_t, t)X^{(1)}_t\right] 
   + \cdots.
   \]
   
   From \eqref{x1} and \eqref{u1}, and writing ${\cal V} = \epsilon^2
   {\cal V}^{(2)} + \cdots$, we find
   \[\label{driftV2}
    {\cal V}^{(2)}
    = \int_0^\infty
    E\left[
      f'(X^{(0)}_t, t)
      \(
     f(X^{(0)}_{t-\beta}, t-\beta) - \lambda U^{(1)}_{t-\beta} 
     \)
    \right]
    d\beta.
   \]
   Again, extending the integral to $\infty$ removes starting
   transients, and the contribution from $(T,\infty)$ is exponentially
   small.  To evaluate this, we need
   \[\label{ff1}
   E\left[
     f'(X^{(0)}_t, t)
     f(X^{(0)}_{t-\beta}, t-\beta)
     \right]
   =
   - \half u^2 k E\left[\sin\(k (X^{(0)}_t - X^{(0)}_{t-\beta}) - \omega \beta\)\right],
   \]
   as the other terms in the complex exponentials average to zero.  As
   $X^{(0)}_t - X^{(0)}_{t-\beta}$ is Gaussian and has zero mean, this
   reduces to
   \[\label{ff2}
   E\left[
     f'(X^{(0)}_t, t)
     f(X^{(0)}_{t-\beta}, t-\beta)
     \right]
   =
   \half u^2 k \sin(\omega\beta) 
   \exp\(-\half k^2 E[(X^{(0)}_t - X^{(0)}_{t-\beta})^2]\),
   \]
   where, from \eqref{XYsol0},
   \[
   E[(X^{(0)}_t - X^{(0)}_{t-\beta})^2]
   = \sigma^2\(\beta + \lambda^{-1}(1 - e^{-\lambda\beta})\).
   \]
   To complete the evaluation of \eqref{driftV2}, we use \eqref{u1} to show
   \[
   \begin{aligned}
   &E\left[
     f'(X^{(0)}_t, t)
     U^{(1)}_{t-\beta}
     \right]
   =
   \half\int_0^\infty
   \exp(-\lambda\alpha) \\
   &\times E\left[f'(X^{(0)}_t, t)
     \(f(X^{(0)}_{t-\alpha-\beta}, t-\alpha-\beta)
     - f(Y^{(0)}_{t-\alpha-\beta}, t-\alpha-\beta)\)\right]
   d\alpha.
   \end{aligned}
   \]
   The first term in the expectation on the right-hand side is as in
   \eqref{ff1} but with different times, so we are essentially left
   with finding
   \[
   K(\alpha+\beta) \equiv E\left[f'(X^{(0)}_t, t)
     f(Y^{(0)}_{t-\alpha-\beta}, t-\alpha-\beta)\right],
   \]
   which itself is only slightly different from \eqref{ff1}.
   As $X^{(0)}_t - Y^{(0)}_{t-\tau}$ is Gaussian with zero mean, we find
   \[\label{ff3}
   K(\tau) = \half u^2 k \sin(\omega\tau)
   \exp\(-\half k^2 E[(X^{(0)}_t - Y^{(0)}_{t-\tau})^2]\),
   \]
   where $\tau \equiv \alpha + \beta$. From \eqref{XYsol0}, 
   \[
   E[(X^{(0)}_t - Y^{(0)}_{t-\tau})^2]
   =
   \sigma^2 \(\tau + \lambda^{-1} (1 + e^{-\lambda\tau})\).
   \]

   We are now ready to evaluate 
   \[
   M(\beta) 
   \equiv 
   E\left[
     f'(X^{(0)}_t, t)
     \(
     f(X^{(0)}_{t-\beta}, t-\beta) - \lambda U^{(1)}_{t-\beta}
     \)
     \right]
   \]
   from \eqref{driftV2}. Gathering together what we have just
   evaluated in \eqref{ff2} and \eqref{ff3}, we find
   \[\label{M}
   \begin{aligned}
   \frac{M(\beta)}{\half u^2 k} &= \sin(\omega\beta)
   \exp\(-\half k^2 \sigma^2\(\beta +
   \lambda^{-1}\(1-e^{-\lambda\beta}\)\)\) \\
   &- \lambda \int_0^\infty
   \exp\(-\lambda\alpha 
   - \half k^2 \sigma^2 
   \(\alpha + \beta + \lambda^{-1}\)\)
   \sin(\omega(\alpha+\beta))
   \sinh\(\frac{k^2 \sigma^2}{2\lambda} e^{-\lambda(\alpha+\beta)}\) d\alpha.
   \end{aligned}
   \]
   Thus \eqref{drift3} and \eqref{driftV2} imply that the Stokes'
   drift to leading order in $\epsilon$ is 
   \[\label{thedrift}
    {\cal V}
    = \epsilon^2 \int_0^\infty M(\beta) d\beta + O(\epsilon^4).
   \]
   Note that in the case of sinusoidal forcing the first correction is
   $O(\epsilon^4)$ rather than $O(\epsilon^3)$ as replacing $\epsilon$
   by $-\epsilon$ is the same as a phase shift, so does not change the
   Stokes' drift.  In fact, all terms with odd powers of $\epsilon$
   are zero, in this case.

   \section{Strong and weak spring limits}
   \label{springlimits}
   
   We briefly consider the large and small $\lambda$ limits.  In the
   large $\lambda$ limit, the particles coalesce, and we expect to
   recover the point particle result\cite{Jansons}. As
   $\lambda/\omega\to\infty$, with $k^2\sigma^2/\omega$ fixed,
   \[
   M \to \half u^2 k \sin(\omega\beta)\exp(-\half k^2 \sigma^2 \beta),
   \]
   which reduces \eqref{thedrift} to
   \[
   {\cal V} = \epsilon^2 \frac{2 u^2 k \omega}{k^4 \sigma^4 + 4\omega^2} + \cdots.
   \]

   On the other hand, small $\lambda$ means that the spring is weak
   and the particles can be separated by many wavelengths. If the
   dumbbell is being used to model a flexible polymer, this is not a
   very physical limit, but we mention it here as this dumbbell model
   may have applications outside of fluid mechanics\cite{Jansons}.

   In the limit $\lambda/\omega\to0$, with $k^2\sigma^2/\omega$ fixed,
   we find
   \[
   M \to \half u^2 k \sin(\omega\beta)\exp(-k^2 \sigma^2 \beta),
   \]
   which reduces \eqref{thedrift} to
   \[
   {\cal V} = \half \epsilon^2 \frac{u^2 k \omega}{k^4 \sigma^4 +
     \omega^2} + \cdots.
   \]
   Note that this has the same form as the single particle result but
   with different constants.

   \section{Numerical results}
   \label{results}

   We already know that we can separate particles with different
   diffusivities using stochastic Stokes' drift\cite{Jansons}, and we
   could always choose wave parameters to make the flexibility of the
   particle unimportant.  So we present here numerical results for
   particles that all have the same diffusivity to highlight the
   scope for sorting particles where only the flexibilities are different.
   
\begin{figure}
     \centering
     \includegraphics[angle=-90,width=\figurewidth]{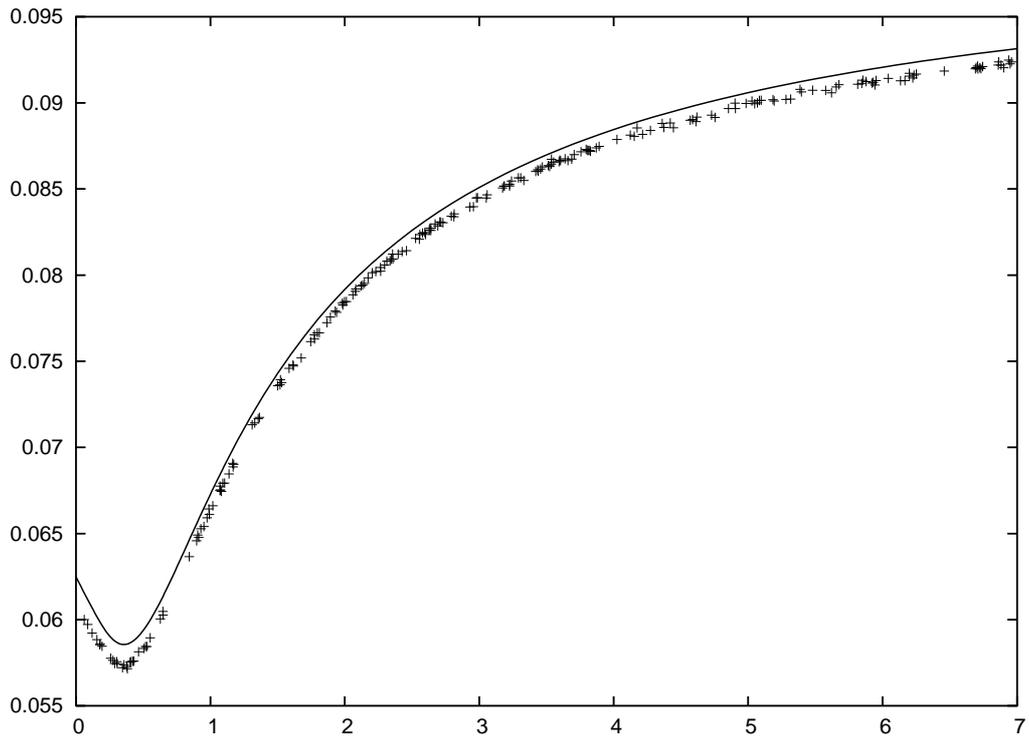}
     \caption{Stokes' drift against $\lambda$ for $\epsilon=0.5$ and
       all other parameters equal to $1$. The curve is the
       leading-order asymptotic result and the data points are from a
       Monte Carlo experiment, with time steps of $0.001$ and run to a
       time $t = 10^8$.}
     \label{fig1}
   \end{figure}

   A graph of the leading-order asymptotic result for the Stokes'
   drift as a function of the spring strength ($\lambda$), with all
   other parameters set to $1$, is compared to a Monte Carlo
   simulation using a stochastic Euler integration with time steps of
   $0.001$ and taken up to $t=10^8$ (see figure \ref{fig1}). The
   agreement is good even when the `small' expansion parameter
   $\epsilon=0.5$. Other parameter values were also tested in the same
   way with equally good agreement.

   Note that in sorting applications, the characteristics of the wave
   motion would be, in general, under user control.
   
   \section{Conclusions}

   We have shown that it is possible to have a Stokes' drift that is
   strongly dependent on the spring constant of the dumbbell particles
   even when all the particle diffusivities are the same.  Thus, in
   principle, it is possible to use this effect to sort particles
   solely by their flexibilities.
   
   The leading-order asymptotic result for the Stokes' drift of an
   elastic dumbbell agrees well with simulations even for quite large
   values of the `small' expansion parameter.  Furthermore, the
   dependence of the Stokes' drift is \emph{not} necessarily monotonic
   in the spring constant $\lambda$, as shown in figure
   \ref{fig1}. So, if we replace $f$ by $f = u\cos(k x - \omega t
   +\phi) - u_0$, there is a range of $u_0$ for which, for sufficiently
   small or sufficiently large $\lambda$, the dumbbell drifts with the
   wave, but for intermediate $\lambda$, the dumbbell drifts against
   the wave.  This begs the question of whether there is some
   mathematically similar flexibility effect in other applications,
   for example in molecular motors.
 
   For the same reasons as discussed in \cite{Jansons}, these
   one-dimensional results can be interpreted as components of
   higher-dimensional results, as at order $\epsilon^2$, cross-terms
   for waves of either different spatial or temporal frequencies
   average to zero.  Thus the `fanout' of particles of different
   diffusivities discussed in \cite{Jansons} under the action of
   several waves would extend to dumbbells of the same diffusivities
   but different spring constants.

   It remains to be seen whether this effect would be significant
   enough to have practical applications.  Note, however, the Brownian
   forcing does not necessarily have to be real thermal agitation, but
   could be some other random forcing, for example eddy diffusion, in
   which case the dumbbell would represent a much bigger flexible
   body.

\newcommand\REVIEW[4]{\textit{#1}, \textbf{#2}, #3, #4}
\newcommand\Name{}
\newcommand\Year{}

   \begin{verbatim}
     $Id: DumbbellDrift.tex 3371 2007-03-22 09:25:08Z kalvis $
   \end{verbatim}

\end{document}